\documentclass[reqno]{amsart}

\usepackage[latin1]{inputenc}
\usepackage{amssymb}
\usepackage{graphicx}
\usepackage{amscd}
\usepackage[hidelinks]{hyperref}
\usepackage{color}
\usepackage{float}
\usepackage{graphics,amsmath,amssymb}
\usepackage{amsthm}
\usepackage{amsfonts}
\usepackage{latexsym}
\usepackage{epsf}
\usepackage{enumerate}
\usepackage{xifthen}
\usepackage{mathrsfs}
\usepackage{dsfont}
\usepackage{makecell}
\usepackage[FIGTOPCAP]{subfigure}
\usepackage{amsmath}
\allowdisplaybreaks[4]
\usepackage{listings}
\usepackage{etoolbox}
\usepackage{fancyhdr}

 \newtheoremstyle{mytheorem}
 {3pt}
 {3pt}
 {\slshape}
 {}
 {\bfseries}
 {.}
 { }
 {}

\numberwithin{equation}{section}

\theoremstyle{theorem}
\newtheorem{theorem}{Theorem}[section]
\newtheorem{corollary}[theorem]{Corollary}
\newtheorem{lemma}[theorem]{Lemma}

\theoremstyle{definition}

\newtheorem{remark}{Remark}[section]

\newcommand{\Keywords}[1]{\ifthenelse{\isempty{#1}}{}{\smallskip \smallskip \noindent \textbf{Keywords}. #1}}
\newcommand{\MSC}[2][2010]{\ifthenelse{\isempty{#2}}{}{\smallskip \smallskip \noindent \textbf{#1MSC}. #2}}
\newcommand{\abstractnote}[1]{\ifthenelse{\isempty{#1}}{}{\smallskip \smallskip \noindent \textsuperscript{\dag}#1}}

\makeatletter
\def\specialsection{\@startsection{section}{1}%
  \z@{\linespacing\@plus\linespacing}{.5\linespacing}%
  {\normalfont}}
\def\section{\@startsection{section}{1}%
  \z@{.7\linespacing\@plus\linespacing}{.5\linespacing}%
  {\normalfont\scshape}}
\patchcmd{\@settitle}{\uppercasenonmath\@title}{\Large\boldmath}{}{}
\patchcmd{\@settitle}{\begin{center}}{\begin{flushleft}}{}{}
\patchcmd{\@settitle}{\end{center}}{\end{flushleft}}{}{}
\patchcmd{\@setauthors}{\MakeUppercase}{\normalsize}{}{}
\patchcmd{\@setauthors}{\centering}{\raggedright}{}{}
\patchcmd{\section}{\scshape}{\large\bfseries\boldmath}{}{}
\patchcmd{\subsection}{\bfseries}{\bfseries\boldmath}{}{}
\renewcommand{\@secnumfont}{\bfseries}
\patchcmd{\@startsection}{\@afterindenttrue}{\@afterindentfalse}{}{}
\patchcmd{\abstract}{\leftmargin3pc}{\leftmargin1pc}{}{}

\def\maketitle{\par
  \@topnum\z@ 
  \@setcopyright
  \thispagestyle{empty}
  \ifx\@empty\shortauthors \let\shortauthors\shorttitle
  \else \andify\shortauthors
  \fi
  \@maketitle@hook
  \begingroup
  \@maketitle
  \toks@\@xp{\shortauthors}\@temptokena\@xp{\shorttitle}%
  \toks4{\def\\{ \ignorespaces}}
  \edef\@tempa{%
    \@nx\markboth{\the\toks4
      \@nx\MakeUppercase{\the\toks@}}{\the\@temptokena}}%
  \@tempa
  \endgroup
  \c@footnote\z@
  \@cleartopmattertags
}
\makeatother

\newcommand{\op}{\overline{p}}
\newcommand{\p}{\overline{P}}
\newcommand{\qbinom}[2]{\begin{bmatrix}#1\\#2\end{bmatrix}}
\newcommand{\oqbinom}[2]{\overline{\begin{bmatrix}#1\\#2\end{bmatrix}}}

\title[Overpartitions with bounded differences]{An overpartition analogue of partitions with bounded differences between largest and smallest parts}

\author[S. Chern]{Shane Chern}
\address{Department of Mathematics, Pennsylvania State University, University Park, PA 16802, USA}
\email{shanechern@psu.edu; chenxiaohang92@gmail.com}

\date{}

\begin{document}

%

\maketitle

\begin{abstract}

We study the generating function for overpartitions with bounded differences between largest and smallest parts, which is analogous to a result of Breuer and Kronholm on integer partitions. We also connect this problem with over $q$-binomial coefficients introduced by Dousse and Kim.

\Keywords{Overpartition, bounded difference between largest and smallest parts, over $q$-binomial coefficient.}

\MSC{Primary 05A17; Secondary 11P84.}
\end{abstract}

\section{Introduction}

A \textit{partition} of a positive integer $n$ is a non-increasing sequence of positive integers $(\lambda_1,\lambda_2,\ldots,\lambda_k)$ with $n=\lambda_1+\lambda_2+\ldots+\lambda_k$. For example, $4$ has five partitions: $4$, $3+1$, $2+2$, $2+1+1$, $1+1+1+1$. Let $p(n)$ denote the number of partitions of $n$. It is well known that its generating function is
$$1+\sum_{n\ge 1}p(n)q^n=\frac{1}{(q)_\infty},$$
where, as usual, we use the standard notations
$$(a)_n=(a;q)_n:=\prod_{k=0}^{n-1} (1-a q^k),$$
and
$$(a)_\infty=(a;q)_\infty:=\lim_{n \to \infty}(a)_n.$$

Recently, Andrews, Beck, and Robbins \cite{ABR2015} considered partitions where the difference between largest and smallest parts is a fixed integer $t$. Let $p(n,t)$ be the number of such partitions of $n$. We have, for example, $p(4,1)=1$ since $4$ has only one such partition: $2+1+1$. In fact, Andrews \textit{et al.} showed that $p(n,0)=d(n)$ and $p(n,1)=n-d(n)$ where, and in the sequel, $d(n)$ denotes the number of divisors of $n$. For $t\ge 2$, they obtained the following generating function
\begin{equation}
\sum_{n\ge 1}p(n,t) q^n= \frac{q^{t-1}(1-q)}{(1-q^t)(1-q^{t-1})}-\frac{q^{t-1}(1-q)}{(1-q^t)(1-q^{t-1})(q)_t}+\frac{q^t}{(1-q^{t-1})(q)_t}.
\end{equation}

In a subsequent paper \cite{BK2016}, Breuer and Kronholm further considered partitions where the difference between largest and smallest parts is at most $t$. When $t=1$, the four such partitions of $4$ are: $4$, $2+2$, $2+1+1$, $1+1+1+1$. Let $p_t(n)$ be the number of such partitions of $n$. Its generating function is, in fact, even neater.  Breuer and Kronholm proved in two ways (one is geometric and the other is combinatorial) that for $t\ge 1$,
\begin{equation}\label{eq:BK}
\sum_{n\ge 1}p_t(n) q^n=\frac{1}{1-q^t}\left(\frac{1}{(q)_t}-1\right).
\end{equation}
Later on, Chapman \cite{Cha2016} also provided an analytic proof which only uses elementary $q$-series manipulation as deep as the $q$-binomial theorem.

In this paper, we shall study an overpartition analogue of Breuer and Kronholm's result.

\section{An overpartition analogue}

An \textit{overpartition} of $n$ is a partition of $n$ where the first occurrence of each distinct part may be overlined. For example, $4$ has fourteen overpartitions:
$$4,\ \overline{4},\ 3+1,\ \overline{3}+1,\ 3+\overline{1},\ \overline{3}+\overline{1},\ 2+2,\ \overline{2}+2,$$
$$2+1+1,\ \overline{2}+1+1,\ 2+\overline{1}+1,\ \overline{2}+\overline{1}+1,\ 1+1+1+1,\ \overline{1}+1+1+1.$$
Overpartitions have many applications in combinatorics \cite{CL2004}, mathematical physics \cite{FJM2005}, representation theory \cite{KK2004}, and symmetric functions \cite{Ber2006}. We denote by $\op(n)$ the number of overpartitions of $n$. It is known that
$$1+\sum_{n\ge 1}\op(n) q^n= \frac{(-q)_\infty}{(q)_\infty}.$$

Given a non-negative integer $t$, let $\op_t(n)$ denote the number of overpartitions of $n$ whose difference between largest and smallest parts is at most $t$. Furthermore, let $g_t(n)$ be the number of such overpartitions of $n$ with one more restriction: if the difference between largest and smallest parts is exactly $t$, then the largest parts cannot be overlined. For example, we have $\op_1(4)=10$:
$$4,\ \overline{4},\ 2+2,\ \overline{2}+2,\ 2+1+1,\ \overline{2}+1+1,$$
$$2+\overline{1}+1,\ \overline{2}+\overline{1}+1,\ 1+1+1+1,\ \overline{1}+1+1+1,$$
and $g_1(4)=8$:
$$4,\ \overline{4},\ 2+2,\ \overline{2}+2,\ 2+1+1,$$
$$2+\overline{1}+1,\ 1+1+1+1,\ \overline{1}+1+1+1.$$
It is clear that $\op_0(n)=2d(n)$ and $g_0(n)=d(n)$. We now write
$$\p_t(q):=\sum_{n\ge 1}\op_t(n) q^n\quad \text{and}\quad G_t(q):=\sum_{n\ge 1}g_t(n) q^n.$$

\begin{theorem}\label{th:1}
For $t\ge 1$, we have
\begin{equation}\label{eq:th1}
G_t(q)=\frac{1}{1-q^t}\left(\frac{(-q)_t}{(q)_t}-1\right).
\end{equation}
\end{theorem}

\begin{remark}
It is of interest to compare the similarity between this generating function and the result of Breuer and Kronholm (Eq. \eqref{eq:BK}).
\end{remark}

\begin{theorem}\label{th:2}
For $t\ge 0$, we have
\begin{equation}\label{eq:th2}
\p_t(q)=2(-1)^t \left(\sum_{n\ge 1}d(n)q^n + \sum_{n=1}^t \frac{(-1)^n}{1-q^n}\left(\frac{(-q)_n}{(q)_n}-1\right)\right).
\end{equation}
\end{theorem}

Before presenting our proofs, we need several $q$-series identities. At first, we introduce the ${}_r\phi_s$ functions:
$${}_{r+1}\phi_s\left(\begin{matrix} a_0,a_1,a_2\ldots,a_r\\ b_1,b_2,\ldots,b_s \end{matrix}; q, z\right):=\sum_{n\ge 0}\frac{(a_0;q)_n(a_1;q)_n\cdots(a_r;q)_n}{(q;q)_n(b_1;q)_n\cdots (b_s;q)_n}\left((-1)^n q^{\binom{n}{2}}\right)^{s-r}z^n.$$

\begin{lemma}[First $q$-Chu--Vandermonde Sum {\cite[Eq. (17.6.2)]{And2010}}]\label{le:chu}
We have
\begin{equation}\label{eq:chu}
{}_{2}\phi_{1}\left(\begin{matrix} a,q^{-n}\\ c \end{matrix}; q, cq^{n}/a\right)=\frac{(c/a;q)_n}{(c;q)_n}.
\end{equation}
\end{lemma}

\begin{lemma}[{\cite[Eq. (17.9.6)]{And2010}}]\label{le:32}
We have
\begin{equation}\label{eq:32}
{}_{3}\phi_{2}\left(\begin{matrix} a,b,c\\ d,e \end{matrix}; q, de/(abc)\right)=\frac{(e/a;q)_n(de/(bc);q)_n}{(e;q)_n(de/(abc);q)_n}{}_{3}\phi_{2}\left(\begin{matrix} a, d/b, d/c\\ d, de/(bc) \end{matrix}; q, e/a\right).
\end{equation}
\end{lemma}

\begin{proof}[Proof of Theorem \ref{th:1}]
Fix the smallest part to be $m\ge 1$. It is easy to see that, for $t\ge 1$, the generating function for overpartitions, where the smallest part is $m$, the largest part is at most $m+t$, and if the largest part is exactly $m+t$, then the largest parts cannot be overlined, is
$$\frac{2q^m}{1-q^m}\frac{1+q^{m+1}}{1-q^{m+1}}\cdots\frac{1+q^{m+t-1}}{1-q^{m+t-1}}\frac{1}{1-q^{m+t}}.$$
Hence,
\begin{align*}
\frac{G_t(q)}{2}&=\sum_{m\ge 1}\frac{q^m}{1-q^m}\frac{1+q^{m+1}}{1-q^{m+1}}\cdots\frac{1+q^{m+t-1}}{1-q^{m+t-1}}\frac{1}{1-q^{m+t}}\\
&=\sum_{m\ge 1}\frac{(q)_{m-1}(-q)_{m+t-1}}{(q)_{m+t}(-q)_m}q^m\\
& =q\sum_{m\ge 0}\frac{(q)_{m}(-q)_{m+t}}{(q)_{m+t+1}(-q)_{m+1}}q^m\\
&=\frac{q(-q)_{t}}{(1+q)(q)_{t+1}}\sum_{m\ge 0}\frac{(q)_{m}(q)_{m}(-q^{t+1})_{m}}{(q)_m(q^{t+2})_{m}(-q^2)_{m}}q^m\\
&=\frac{q(-q)_{t}}{(1+q)(q)_{t+1}}\ {}_{3}\phi_{2}\left(\begin{matrix} q,q,-q^{t+1}\\ -q^2,q^{t+2} \end{matrix}; q, q\right)\\
\text{(by Eq. \eqref{eq:32})}& = \frac{q(-q)_{t}}{(1+q)(q)_{t+1}}\frac{(q^{t+1})_\infty(q^2)_\infty}{(q^{t+2})_\infty(q)_\infty}\ {}_{3}\phi_{2}\left(\begin{matrix} q,-q,q^{1-t}\\ -q^2,q^{2} \end{matrix}; q, q^{t+1}\right)\\
&=  \frac{q(-q)_{t}}{(1-q^2)(q)_{t}} \sum_{m\ge 0}\frac{(-q)_{m}(q^{1-t})_{m}}{(-q^2)_{m}(q^2)_m}q^{m(t+1)}\\
&=  \frac{q(-q)_{t}}{(1-q^2)(q)_{t}} \sum_{m\ge 0}\frac{\frac{(-1)_{m+1}}{2}\frac{(q^{-t})_{m+1}}{1-q^{-t}}}{\frac{(-q)_{m+1}}{1+q}\frac{(q)_{m+1}}{1-q}}q^{(m+1)(t+1)}q^{-(t+1)}\\
&=  -\frac{(-q)_{t}}{2(1-q^t)(q)_{t}} \left({}_{2}\phi_{1}\left(\begin{matrix} -1,q^{-t}\\ -q \end{matrix}; q, q^{t+1}\right)-1\right)\\
\text{(by Eq. \eqref{eq:chu})}&=  -\frac{(-q)_{t}}{2(1-q^t)(q)_{t}} \left(\frac{(q)_t}{(-q)_t}-1\right)\\
&= \frac{1}{2(1-q^t)}\left(\frac{(-q)_t}{(q)_t}-1\right).
\end{align*}
We therefore have
$$G_t(q)=\frac{1}{1-q^t}\left(\frac{(-q)_t}{(q)_t}-1\right).$$
This ends the proof of Theorem \ref{th:1}.
\end{proof}

\begin{proof}[Proof of Theorem \ref{th:2}]
We first note that
$$\p_0(q)=2\sum_{n\ge 1}d(n)q^n,$$
and for $t\ge 0$,
$$\p_t(q)=2\sum_{m\ge 1}\frac{q^m}{1-q^m}\frac{1+q^{m+1}}{1-q^{m+1}}\cdots\frac{1+q^{m+t}}{1-q^{m+t}}.$$
We also have, for $n\ge 1$,
\begin{align*}
\p_n(q)+\p_{n-1}(q)&=2\sum_{m\ge 1}\frac{q^m}{1-q^m}\frac{1+q^{m+1}}{1-q^{m+1}}\cdots\frac{1+q^{m+n-1}}{1-q^{m+n-1}}\left(\frac{1+q^{m+n}}{1-q^{m+n}}+1\right)\\
&=4\sum_{m\ge 1}\frac{q^m}{1-q^m}\frac{1+q^{m+1}}{1-q^{m+1}}\cdots\frac{1+q^{m+n-1}}{1-q^{m+n-1}}\frac{1}{1-q^{m+n}}\\
&=2G_n(q).
\end{align*}
Hence, for any $t\ge 0$, it follows that
\begin{align*}
\p_t(q)&=(-1)^t\left(\p_0(q)+2\sum_{n=1}^t (-1)^n G_n(q)\right)\\
&=2(-1)^t \left(\sum_{n\ge 1}d(n)q^n + \sum_{n=1}^t \frac{(-1)^n}{1-q^n}\left(\frac{(-q)_n}{(q)_n}-1\right)\right).
\end{align*}
\end{proof}

\begin{corollary}
For any $t\ge 0$ and $n\ge 1$, $\op_t(n)$ is an even integer. Furthermore, $\op_t(n)$ is divisible by $4$ if and only if $n$ is not a perfect square.
\end{corollary}

\begin{proof}
From the previous proof, we see that
\begin{align*}
\p_t(q)&=(-1)^t\left(\p_0(q)+2\sum_{n=1}^t (-1)^n G_n(q)\right)\\
&=(-1)^t \Bigg(2\sum_{n\ge 1}d(n)q^n\\
&\quad\quad\quad\quad\quad+4\sum_{n=1}^t (-1)^n \sum_{m\ge 1}\frac{q^m}{1-q^m}\frac{1+q^{m+1}}{1-q^{m+1}}\cdots\frac{1+q^{m+n-1}}{1-q^{m+n-1}}\frac{1}{1-q^{m+n}}\Bigg).
\end{align*}
It follows that
$$\op_t(n)\equiv 2d(n) \pmod{4}.$$
We finally note that $d(n)$ is odd if and only if $n$ is a perfect square. This completes the proof.
\end{proof}

\section{The viewpoint of over $q$-binomial coefficients}

The \textit{$q$-binomial coefficient}, also known as \textit{Gaussian polynomial}, is defined as
$$\qbinom{M+N}{N}=\qbinom{M+N}{N}_q:=\frac{(q)_{M+N}}{(q)_{M}(q)_{N}}.$$
We know that it is the generating function for partitions where the largest part is at most $M$ and the number of parts is at most $N$. In a recent paper \cite{DK2017}, Dousse and Kim introduced the \textit{over $q$-binomial coefficient}, denoted by
$$\oqbinom{M+N}{N}=\oqbinom{M+N}{N}_q,$$
which is an overpartition analogue of $q$-binomial coefficient defined as the generating function for overpartitions where the largest part is at most $M$ and the number of parts is at most $N$. They showed that for positive integers $M$ and $N$
$$\oqbinom{M+N}{N}=\sum_{k=0}^{\min(M,N)}q^{\binom{k+1}{2}}\frac{(q)_{M+N-k}}{(q)_k(q)_{M-k}(q)_{N-k}}.$$
Of course, if we agree that the number of such overpartitions of $0$ is one, then this identity also holds for $M=0$ or $N=0$.

Over $q$-binomial coefficients have many properties similar to those of the standard $q$-binomial coefficients. For example, the following recurrence relation
\begin{equation}\label{eq:ob}
\oqbinom{M+N}{N}=\oqbinom{M+N-1}{N-1}+q^N \oqbinom{M+N-1}{N}+ q^N \oqbinom{M+N-2}{N-1}
\end{equation}
holds for any positive integers $M$ and $N$ (see \cite[Eq. (1.1)]{DK2017}). In fact, it can be proved combinatorially.

Motivated by \cite{Cha2016}, we provide an alternative proof of Theorem \ref{th:1} using over $q$-binomial coefficients. It avoids the application of complicated $q$-series identities such as the $q$-Chu--Vandermonde sum.

\begin{proof}[Second proof of Theorem \ref{th:1}]
Here we always assume $t$ to be a positive integer. Let $\lambda=(\lambda_1,\ldots,\lambda_r)$ be an overpartition of $n$ with exactly $r$ parts, $\lambda_r=m\ge 1$, and $\lambda_1\le m+t$. Then $\mu=(\lambda_1 -m,\ldots,\lambda_{r-1}-m)$ is an overpartition of $n-rm$ with at most $r-1$ parts and greatest part $\le t$. Note that the first occurrence of the smallest part of $\lambda$ can be either overlined or not. Hence the generating function for such overpartitions is
$$2q^{rm}\oqbinom{t+r-1}{t},$$
and hence
$$\p_t(q)=2\sum_{r\ge 1}\sum_{m\ge 1}q^{rm}\oqbinom{t+r-1}{t}=2\sum_{r\ge 1}\frac{q^r}{1-q^r}\oqbinom{t+r-1}{t}.$$
We remark that this identity also holds for $t=0$.

On the other hand, overpartitions where the difference between largest and smallest parts is at most $t$ can be divided into three disjoint cases:
\begin{enumerate}
\item The largest part is at most $t$;
\item The largest part is greater than $t$, the difference between largest and smallest parts is exactly $t$, and the first occurrence of the smallest part is overlined;
\item Otherwise.
\end{enumerate}
For Case (1), one readily sees the generating function is
$$\frac{(-q)_t}{(q)_t}-1.$$
For Case (2), its generating function is
$$2\sum_{m\ge 1} \frac{q^m}{1-q^m}\frac{1+q^{m+1}}{1-q^{m+1}}\cdots\frac{1+q^{m+t-1}}{1-q^{m+t-1}}\frac{q^{m+t}}{1-q^{m+t}}=\frac{\p_t(q)-\p_{t-1}(q)}{2}.$$
Finally, let $\lambda=(\lambda_1,\ldots,\lambda_r)$ be an overpartition of Case (3) with $\lambda_1=m+t$ (and so $m\ge 1$). We note that $\mu=(\lambda_1 -m,\ldots,\lambda_{r}-m)$ is an overpartition of $|\lambda|-rm$ with at most $r$ parts and largest part being exactly $t$. Hence the generating function is 
$$\sum_{r\ge 1}\sum_{m\ge 1}q^{rm}\left(\oqbinom{t+r}{t}-\oqbinom{t+r-1}{t-1}\right)=\sum_{r\ge 1}\frac{q^r}{1-q^r}\left(\oqbinom{t+r}{t}-\oqbinom{t+r-1}{t-1}\right).$$
We therefore have
\begin{align*}
\p_t(q)&=\left(\frac{(-q)_t}{(q)_t}-1\right)+\frac{\p_t(q)-\p_{t-1}(q)}{2}\\
&\quad+\sum_{r\ge 1}\frac{q^r}{1-q^r}\left(\oqbinom{t+r}{t}-\oqbinom{t+r-1}{t-1}\right).
\end{align*}

Now we take $M\to r$ and $N\to t$ in Eq. \eqref{eq:ob} and rewrite it as
$$\oqbinom{t+r}{t}-\oqbinom{t+r-1}{t-1}=q^t \left(\oqbinom{t+r-1}{t}+ \oqbinom{t+r-2}{t-1}\right).$$
We then multiply both sides by $q^r/(1-q^r)$ and sum over $r$
\begin{align*}
&\sum_{r\ge 1}\frac{q^r}{1-q^r}\left(\oqbinom{t+r}{t}-\oqbinom{t+r-1}{t-1}\right)\\
&\quad =q^t\left(\sum_{r\ge 1}\frac{q^r}{1-q^r}\oqbinom{t+r-1}{t}+ \sum_{r\ge 1}\frac{q^r}{1-q^r}\oqbinom{t+r-2}{t-1}\right).
\end{align*}
From the foregoing argument, we therefore have
\begin{align*}
\p_t(q)-\left(\frac{(-q)_t}{(q)_t}-1\right)-\frac{\p_t(q)-\p_{t-1}(q)}{2}=q^t \frac{\p_t(q)+\p_{t-1}(q)}{2}.
\end{align*}
Hence
$$G_t(q)=\frac{\p_t(q)+\p_{t-1}(q)}{2}=\frac{1}{1-q^t}\left(\frac{(-q)_t}{(q)_t}-1\right).$$
\end{proof}

\subsection*{Acknowledgements}

I would like to thank George E. Andrews and Ae Ja Yee for some helpful discussions.

\bibliographystyle{amsplain}

\end{document}